\documentstyle[12pt]{article}

\newlength{\abstractwidth}
\renewcommand{\thefootnote}{\fnsymbol{footnote}}
\renewcommand{\thanks}[1]{\footnote{#1}} 
\newcommand{\starttext}{ \setcounter{footnote}{0}
\renewcommand{\thefootnote}{\arabic{footnote}}}

\newcommand{\be}{\begin{equation}}
\newcommand{\bea}{\begin{eqnarray}}
\newcommand{\eea}{\end{eqnarray}} 
\newcommand{\bean}{\begin{eqnarray*}}
\newcommand{\eean}{\end{eqnarray*}} 
\newcommand{\ee}{\end{equation}}
 \newcommand{\<}{\langle}
\renewcommand{\>}{\rangle}
\def\ba{\begin{eqnarray}}
\def\ea{\end{eqnarray}}

\def\o{\omega}
\def\Re{{\rm Re}}

\def\log{\,{\rm log}\,}

\def\o{\omega}

\def\e{\varepsilon}

\def\o{\omega}

\def\p{\partial}

\def\R{{\bf R}}

\def\ddb{{\partial\bar\partial}}

\def\F{{\cal F}}

\def\[{{\bf [}}
\def\]{{\bf ]}}

\begin{document}
\starttext \baselineskip=15pt \setcounter{footnote}{0}
\newtheorem{theorem}{Theorem}
\newtheorem{lemma}{Lemma}
\newtheorem{definition}{Definition}
\newtheorem{proposition}{Proposition}
\newtheorem{corollary}{Corollary}
\newtheorem{remark}{Remark}

\begin{center}
{\Large \bf A SECOND ORDER ESTIMATE FOR GENERAL COMPLEX HESSIAN EQUATIONS
\footnote{Work supported in part by the National Science Foundation under Grant DMS-1266033, DMS-1605968 and DMS-1308136.
}}

\bigskip

{\large Duong H. Phong, Sebastien Picard and Xiangwen Zhang} \\
\medskip
\begin{abstract}
We derive a priori $C^2$ estimates for the $\chi$-plurisubharmonic solutions of general complex Hessian equations with right-hand side depending on gradients. 
\end{abstract}
\end{center}

\section{Introduction}
\par
Let $(X, \omega)$ be a compact K\"ahler manifold of dimension $n\geq 2$. Let $u\in C^{\infty}(X)$ and consider a $(1, 1)$ form $\chi(z, u)$ possibly depending on $u$ and satisfying the positivity condition $\chi \geq \e \omega$
for some $\e>0$. We define
\be
g=\chi(z, u) + i \ddb u,
\ee
and $u$ is called $\chi$-plurisubharmonic if $g>0$ as a $(1, 1)$ form. In this paper, we are concerned with the following complex Hessian equation, for $1\leq k\leq n$,
\be\label{equation}
\left(\chi(z, u) + i \ddb u \right)^k \wedge \omega^{n-k} = \psi(z, Du, u) \, \omega^n,
\ee
where $\psi(z, v, u)\in C^{\infty}\left((T^{1,0}(X))^* \times \R \right)$ is a given strictly positive function. 

\medskip
The complex Hessian equation can be viewed as an intermediate equation between the Laplace equation and the complex Monge-Amp\`ere equation. It encompasses the most natural invariants of the complex Hessian matrix of a real valued function, namely the elementary symmetric polynomials of its eigenvalues. When $k=1$, equation (\ref{equation}) is quasilinear, and the estimates follow from the classical theory of quasilinear PDE. The real counterparts of (\ref{equation}) for $1<k\leq n$, with $\psi$ not depending on the gradient of $u$, have been studied extensively in the literature (see the survey paper \cite{Wang} and more recent related work \cite{Guanbo}), as these equations appear naturally and play very important roles in both classical and conformal geometry. When the right-hand side $\psi$ depends on the gradient of the solution, even the real case has been a long standing problem due to substantial difficulties in obtaining a priori $C^2$ estimates. This problem was recently solved by Guan-Ren-Wang \cite{GRW} for convex solutions of real Hessian equations.

\medskip
In the complex case, the equation (\ref{equation}) with $\psi=\psi(z, u)$ has been extensively studied in recent years, due to its appearance in many geometric problems, including the $J$-flow \cite{SW} and quaternionic geometry \cite{AV}. The related Dirichlet problem for equation (\ref{equation}) on domains in ${\bf C}^n$ has been studied by Li \cite{Li} and Blocki \cite{Blocki1}. The corresponding problem on compact K\"ahler or Hermitian manifolds has also been studied extensively, see, for example, \cite{DK, Hou, KN, LN,  Zhang}. In particular, as a crucial step in the continuity method, $C^2$ estimates for complex Hessian type equations have been studied in various settings, see \cite{HMW, Sun, Gabor, GTW, DZhang}. 

\medskip
However, the equation (\ref{equation}) with $\psi=\psi(z, Du, u)$ has been much less studied. An important case corresponding to $k=n=2$, so that it is actually a Monge-Amp\`ere equation in two dimensions, is central to the solution by Fu and Yau \cite{FY, FY1} of a Strominger system on a toric fibration over a K3 surface. A natural generalization of this case to general dimension $n$ was suggested by Fu and Yau \cite{FY} and can be expressed as
\begin{equation}\label{fyequ}
\left( \left(e^u + f e^{-u}\right)\omega + n \, i \ddb u\right)^2 \wedge \omega^{n-2} = \psi(z, Du, u) \, \omega^n,
\end{equation}
where $\psi(z, v, u)$ is a function on $(T^{1,0}(X))^*\times {\bf R}$ 
with a particular structure, and $(X, \omega)$ is a compact K\"ahler manifold. A priori estimates for this equation were obtained by the authors in \cite{PPZ}.

\medskip
In this paper, motivated by our previous work \cite{PPZ}, we study a priori $C^2$ estimate for the equation (\ref{equation}) with general $\chi(z,u)$ and general right hand side $\psi(z,Du,u)$. Building on the techniques developed by Guan-Ren-Wang in \cite{GRW} (see also \cite{LRW}) for real Hessian equations), we can prove the following theorem.

\begin{theorem}
Let $(X, \omega)$ be a compact K\"ahler manifold of complex dimension $n$. Suppose $u\in C^4(X)$ is a solution of equation (\ref{equation}) with $g = \chi+ i \ddb u >0$ and $\chi(z, u) \geq \e \omega$. Let $0<\psi(z, v, u)\in C^{\infty}((T^{1,0} X)^* \times \R)$. Then we have the following uniform second order derivative estimate
\begin{eqnarray}
 |D\bar D u|_{\omega} \leq C,
\end{eqnarray}
where $C$ is a positive constant depending only on $\e, n, k$, $\sup_X |u|$, $\sup_{X} | Du|$, and
the $C^2$ norm of $\chi$ as a function of $(u, z)$, the infimum of $\psi$, and the $C^2$ norm of $\psi$ as a function of $(z, Du, u)$, all restricted to the ranges in $Du$ and $u$ defined by the uniform upper bounds on $|u|$ and $|Du|$.

\end{theorem}

\smallskip
We remark that the above estimate is stated for $\chi$-plurisubharmonic solutions, that is, $g = \chi+ i \ddb u >0$. Actually, we only need to assume that $g\in \Gamma_{k+1}$ cone (see (\ref{cone}) below for the definition of the Garding cone $\Gamma_k$ and also the discussion in Remark \ref{remarkcone} at the end of the paper). However, a better condition would be $g\in \Gamma_k$, which is the natural cone for ellipticity. In fact, this is still an open problem even for real Hessian equations when $2<k<n$. If $k=2$, Guan-Ren-Wang \cite{GRW} removed the convexity assumption by investigating the structure of the operator. A simpler argument was given recently by Spruck-Xiao \cite{SX}. However, the arguments are not applicable to the complex case due to the difference between the terms $|DDu|^2$ and $|D\bar D u|^2$ in the complex setting. When $k=2$ in the complex setting, $C^2$ estimates for equation (\ref{fyequ}) were obtained in \cite{PPZ} without the plurisubharmonicity assumption, but the techniques rely on the specific right hand side $\psi(z,Du,u)$ studied there.
\par
 We also note that if $k=n$, the condition $g = \chi+ i \ddb u >0$ is the natural assumption for the ellipticity of equation (\ref{equation}). Thus, our result implies the a priori $C^2$ estimate for complex Monge-Amp\`ere equations with right hand side depending on gradients:
\bean
\left(\chi(z, u) + i \ddb u\right)^n = \psi(z, Du, u) \, \omega^n.
\eean
This generalizes the $C^2$ estimate for the equation studied by
Fu and Yau \cite{FY, FY1} mentioned above, which corresponds to $n=2$
and a specific form $\chi(z,u)$ as well as a specific right hand side $\psi(z,Du,u)$. For dimension $n\geq 2$ and $k=n$, the estimate was obtained by Guan-Ma \cite{Guan} using a different method where the structure of the Monge-Amp\`ere operator plays an important role.

\medskip
Compared to the estimates when $\psi= \psi(z, u)$, the dependence on the gradient of $u$ in the equation (\ref{equation}) creates substantial new difficulties. The main obstacle is the appearance of terms such as $|DDu|^2$ and $|D\bar D u|^2$ when one differentiates the equation twice. We adapt the techniques used in \cite{GRW} and \cite{LRW} for real Hessian equations to overcome these difficulties. Furthermore, we also need to handle properly some subtle issues when dealing with the third order terms due to complex conjugacy.

\medskip

\section{Preliminaries}
\setcounter{equation}{0}

\par
Let $\sigma_k$ be the k-th elementary symmetric function, that is, for $1\leq k\leq n$ and $\lambda= (\lambda_1, \cdots, \lambda_n) \in \R^n$,
\begin{eqnarray*}
 \sigma_k(\lambda) = \sum_{1 < i_1 < \cdots < i_k < n} \lambda_{i_1} \lambda_{i_2}\cdots \lambda_{i_n}.
\end{eqnarray*}
Let $\lambda(a_{\bar j i})$ denote the eigenvalues of a Hermitian symmetric matrix $(a_{\bar j i})$ with respect to the background K\"ahler metric $\o$. We define $ \sigma_k(a_{\bar j i}) = \sigma_k \left( \lambda(a_{\bar j i})\right)$.
This definition can be naturally extended to complex manifolds. Denoting $A^{1, 1}(X)$ to be the space of smooth real $(1, 1)$-forms on a compact K\"ahler manifold $(X, \omega)$, we define for any $g \in A^{1, 1}(X)$, 
\begin{eqnarray*}
\sigma_k (g) ={n\choose k}\frac{ g^k \wedge \omega^{n-k}}{\omega^n}.
\end{eqnarray*}
Using the above notation, we can re-write equation (\ref{equation}) as following:
\begin{eqnarray}\label{sigma_eqn}
\sigma_k(g) = \sigma_k ( \chi_{\bar{j} i} + u_{\bar{j} i})= \psi(z, Du, u).
\end{eqnarray}
We will use the notation
\bean
\sigma_k^{p \bar q} = \frac{\partial \sigma_k(g)}{\partial g_{\bar{q} p}}, \ \  \ \ \sigma_k^{p \bar q, r \bar s}= \frac{\partial^2 \sigma_k(g)}{\partial g_{\bar{q} p} \bar\partial g_{\bar{s} r}}.
\eean
The symbol $D$ will indicate the covariant derivative with respect to the given metric $\omega$. All norms and inner products will be with respect to $\o$ unless denoted otherwise. We will denote by $\lambda_1, \lambda_2, \cdots, \lambda_n$ the eigenvalues of $g_{\bar j i} = \chi_{\bar j i} + u_{\bar j i}$ with respect to $\omega$, and use the ordering $\lambda_1\geq \lambda_2\geq \cdots\geq\lambda_n>0$. Our calculations will be carried out at a point $z$ on the manifold $X$, and we shall use coordinates such that at this point $\o = i\sum \delta_{\ell k} \, dz^k \wedge d\bar{z}^\ell$ and $g_{\bar j i}$ is diagonal. We will also use the notation
\be \nonumber
{\mathcal{F}} = \sum_p \sigma_k^{p \bar{p}}.
\ee
Differentiating equation (\ref{sigma_eqn}) yields
\bea
 \label{diff_once} \sigma_k^{p \bar{q}} D_{\bar{j}} g_{\bar{q} p} &=& D_{\bar{j}} \psi.
 \eea
Differentiating the equation a second time gives
 \bea
 &&\label{diff_twice} \sigma_k^{p \bar{q}} D_i D_{\bar{j}} g_{\bar{q} p} + \sigma_k^{p \bar{q},r \bar{s}} D_i g_{\bar{q} p} D_{\bar{j}} g_{\bar{s} r} = D_i D_{\bar{j}}\psi \nonumber\\
&\geq& -C(1 + |DDu|^2 + |D \bar{D} u|^2) + \sum_\ell \psi_{v_\ell} u_{\ell \bar{j} i} + \sum_\ell \psi_{\bar{v}_\ell} u_{\bar{\ell} \bar{j} i}.
\eea
We will denote by $C$ a uniform constant which depends only on $(X, \omega), n, k$, $\| \chi \|_{C^2}, \inf \psi, \| u \|_{C^1}$ and $\| \psi \|_{C^2}$.
We now compute the operator $\sigma_k^{p \bar{q}} D_p D_{\bar{q}}$ acting on $g_{\bar{j} i} = \chi_{\bar{j} i} + u_{\bar{j} i}$. Recalling that $\chi_{\bar{j} i}$ depends on $u$, we estimate
\bea
\sigma_k^{p \bar{q}} D_p D_{\bar{q}}  g_{\bar{j} i} &=&  \sigma_k^{p \bar{q}} D_p D_{\bar{q}} D_i D_{\bar{j}} u + \sigma_k^{p \bar{q}} D_p D_{\bar{q}} \chi_{\bar{j} i} \nonumber\\
&\geq& \sigma_k^{p \bar{q}} D_p D_{\bar{q}} D_i D_{\bar{j}} u - C(1 + \lambda_1) \F.
\eea
Commuting derivatives
\bea
D_p D_{\bar{q}} D_i D_{\bar{j}} u &=& D_i D_{\bar{j}} D_p D_{\bar{q}} u -  R_{\bar{q} i \bar{j}}{}^{\bar{a}} u_{\bar{a} p} + R_{\bar{q} p \bar{j}}{}^{\bar{a}} u_{\bar{a} i} \nonumber\\
&=& D_i D_{\bar{j}} g_{\bar{q} p} - D_i D_{\bar{j}} \chi_{\bar{q} p} -  R_{\bar{q} i \bar{j}}{}^{\bar{a}} u_{\bar{a} p} + R_{\bar{q} p \bar{j}}{}^{\bar{a}} u_{\bar{a} i}.
\eea
Therefore, by (\ref{diff_twice}),
\bea \label{DDDD}
\sigma_k^{p \bar{q}} D_p D_{\bar{q}} g_{\bar{j} i} &\geq&  - \sigma_k^{p \bar{q},r \bar{s}} D_j g_{\bar{q} p} D_{\bar{j}} g_{\bar{s} r}+ \sum \psi_{v_\ell} g_{\bar{j} i \ell} + \sum \psi_{\bar{v_\ell}}  g_{\bar{j} i \bar{\ell}}\nonumber\\
&& -C(1 + |DDu|^2 + |D \bar{D} u|^2+ (1+\lambda_1) \mathcal{F}).
\eea
We next compute the operator $\sigma_k^{p \bar{q}} D_p D_{\bar{q}}$ acting on $|Du|^2$. Introduce the notation
\be
|DDu|^2_{\sigma \omega}= \sigma_k^{p \bar{q}} \o^{m \bar{\ell}}  D_p D_m u D_{\bar{q}} D_{\bar{\ell}} u, \ \  |D \bar{D} u|^2_{\sigma \omega} = \sigma_k^{p \bar{q}} \o^{m \bar{\ell}}  D_p D_{\bar{\ell}} u D_m D_{\bar{q}} u.
\ee
Then
\bea \label{DD|Du|^2}
\sigma_k^{p \bar{q}} |Du|^2_{\bar{q} p} &=&  \sigma_k^{p \bar{q}} (D_p D_{\bar{q}} D_m u D^m u + D_m u D_p D_{\bar{q}} D^m u)+|DDu|^2_{\sigma \o} + |D \bar{D} u|^2_{\sigma \o} \nonumber\\
&=&  \sigma_k^{p \bar{q}} \{ D_m (g_{\bar{q} p} - \chi_{\bar{q} p}) D^m u + D_m u D^m (g_{\bar{q} p} - \chi_{\bar{q} p}) \}+ \sigma_k^{p \bar{q}} R_{\bar{q} p}{}^{m \bar{\ell}} u_{\bar{\ell}}u_m  
 \nonumber\\
&&+|DDu|^2_{\sigma \o} + |D \bar{D} u|^2_{\sigma \o}.
\eea
Using the differentiated equation we obtain
\bea 
\sigma_k^{p \bar{q}} |Du|^2_{\bar{q} p} &\geq&  2 \Re \< Du, D \psi \> - C (1+{\mathcal{F}}) + |DDu|^2_{\sigma \o} +  |D \bar{D} u|^2_{\sigma \o} \nonumber\\
&\geq& 2 \Re \{ \sum_{p,m} (D_p D_m u D_{\bar p }u + D_p u D_{\bar{p}} D_m u) \psi_{v_m} \} -C(1 + {\mathcal{F}}) +|DDu|^2_{\sigma \o} + |D \bar{D} u|^2_{\sigma \o}. \nonumber
\eea
To simplify the expression, we introduce the notation
\be
\langle D|Du|^2, D_{\bar{v}} \psi \rangle = \sum_m  (D_m D_p u D^p u \, \psi_{v_m} + D_p u D_m D^p u \, \psi_{v_m}).
\ee
We obtain
\be \label{DD|Du|^2}
\sigma_k^{p \bar{q}} |Du|^2_{\bar{q} p} \geq 2 \Re \langle D|Du|^2, D_{\bar{v}} \psi \rangle -C(1 + {\mathcal{F}}) +|DDu|^2_{\sigma \o} +  |D \bar{D} u|^2_{\sigma \o}.
\ee
We also compute
\be \label{sigma_DDu}
- \sigma_k^{p \bar{q}} u_{\bar{q} p} = \sigma_k^{p \bar{q}} (\chi_{\bar{q} p} - g_{\bar{q} p}) \geq  \e {\mathcal{F}} - k \psi.
\ee
\medskip

\section{The $C^2$ estimate}
\setcounter{equation}{0}
\par
In this section, we give the proof of the estimate stated in the theorem. When $k=1$, the equation (\ref{equation}) becomes
\bea
\Delta_\o u+{\rm Tr}_\o\chi(z,u)=n\psi(z,Du,u)
\eea
where $\Delta_\o$ and ${\rm Tr}_\o$ are the Laplacian and trace with respect to the background metric $\o$. It follows that $\Delta_\o u$ is bounded, and the desired estimate follows in turn from the positivity of the metric $g$. Henceforth, we assume that $k\geq 2$.
Motivated by the idea from \cite{GRW} for real Hessian equations, we apply the maximum principle to the following test function:
\begin{eqnarray}\label{testfunction}
G = \log P_m + mN |Du|^2 - mMu,
\end{eqnarray}
where $P_m = \sum_j \lambda_j^m$. Here, $m$, $M$ and $N$ are large positive constants to be determined later. We may assume that the maximum of $G$ is achieved at some point $z\in X$. After rotating the coordinates, we may assume that the matrix $g_{\bar{j} i} = \chi_{\bar{j} i}+ u_{\bar{j} i}$ is diagonal.
\par
Recall that if $F(A)= f(\lambda_1, \cdots, \lambda_n)$ is a symmetric function of the eigenvalues of a Hermitian matrix $A=(a_{\bar ji})$, then at a diagonal matrix $A$ with distinct eigenvalues, we have (see \cite{Ball}), 
\begin{eqnarray}
 \label{firstorder} F^{i\bar j} &=& \delta_{ij} f_i,\\
  \label{secondorder} F^{i\bar j, r\bar s} w_{i\bar j k} w_{r\bar s \bar k} &=& \sum f_{ij} w_{i\bar i k} w_{j\bar j \bar k} + \sum_{p\neq q}\frac{f_p - f_q}{\lambda_p-\lambda_q} | w_{p\bar q k}|^2.
\end{eqnarray}
where $F^{i\bar j}={\p F\over \p a_{\bar ji}}$,
$F^{i\bar j,r\bar s}={\p^2 F\over \p a_{\bar ji}\p a_{\bar sr}}$, and $w_{i\bar jk}$ is an arbitrary tensor.
Using these identities to differentiate $G$, we first obtain the critical equation
\be \label{testfunction1}
{D P_m \over P_m} + mN D|Du|^2 - mM Du =0.
\ee
Differentiating $G$ a second time and contracting with $\sigma_k^{p \bar{q}}$ yields
\bea \label{testfunction2}
0 &\geq& {m \over P_m} \bigg\{ \sum_j \lambda_j^{m-1} \sigma_k^{p \bar{p}} D_p D_{\bar{p}} g_{\bar{j} j} \bigg\} - {|D P_m|^2_\sigma \over P_m^2} + mN \sigma_k^{p \bar{p}} |Du|^2_{\bar{p} p} -mM \sigma_k^{p \bar{p}} u_{\bar{p} p}\nonumber\\
&&+ {m \over P_m} \bigg\{ (m-1)\sum_j \lambda_j^{m-2} \sigma_k^{p \bar{p}} |D_p g_{\bar{j} j}|^2 + \sigma_k^{p \bar{p}} \sum_{i \neq j} {\lambda_i^{m-1} - \lambda_j^{m-1} \over \lambda_i - \lambda_j} |D_p g_{\bar{j} i}|^2 \bigg\}.
\eea
Here, we used the notation
$|\eta|_\sigma^2 = \sigma_k^{p \bar{q}} \eta_p \eta_{\bar{q}}$.
Substituting (\ref{DDDD}), (\ref{DD|Du|^2}) and (\ref{sigma_DDu})
\bea \label{mainineq}
0 &\geq& {1 \over P_m} \bigg\{ -C \sum_j \lambda_j^{m-1} (1 + |DDu|^2 + |D \bar{D} u|^2+ (1+\lambda_1) {\mathcal{F}} ) \bigg\} \nonumber\\
&&+  {1 \over P_m} \bigg\{ \sum_j \lambda_j^{m-1}  (-\sigma_k^{p \bar{q},r \bar{s}} D_j g_{\bar{q} p} D_{\bar{j}} g_{\bar{s} r} + \sum_\ell \psi_{v_\ell} g_{\bar{j} j \ell} + \sum_\ell \psi_{\bar{v}_\ell} g_{\bar{j} j \bar{\ell}}) \bigg\} \nonumber\\
&&+ {1 \over P_m} \bigg\{ (m-1)\sum_j \lambda_j^{m-2} \sigma_k^{p \bar{p}} |D_p g_{\bar{j} j}|^2 + \sigma_k^{p \bar{p}} \sum_{i \neq j} {\lambda_i^{m-1} - \lambda_j^{m-1} \over \lambda_i - \lambda_j} |D_p g_{\bar{j} i}|^2 \bigg\} \nonumber\\
&&- {|D P_m|^2_\sigma \over m P_m^2} + N (|DDu|^2_{\sigma \o} + |D \bar{D} u|^2_{\sigma \o}) \nonumber\\
&&+ N \langle D|Du|^2, D_{\bar{v}} \psi \rangle +  N \langle D_{\bar{v}} \psi , D|Du|^2 \rangle  + (M\e- C N) {\mathcal{F}} - k M \psi.
\eea
From the critical equation (\ref{testfunction1}), we have
\bean
{1 \over P_m} \sum_{j,\ell } \lambda_j^{m-1} g_{\bar{j} j \ell} \psi_{v_\ell} &=& {1 \over m} \langle {D P_m \over P_m},  D_{\bar{v}} \psi \rangle = -N \langle D|Du|^2, D_{\bar{v}} \psi \rangle +M\langle Du, D_{\bar{v}} \psi \rangle.
\eean
It follows that
\bean
&&{1\over P_m}\sum_{j,\ell}\left( \psi_{v_\ell}  g_{\bar{j} j \ell} + \psi_{\bar{v}_\ell} g_{\bar{j} j \bar{\ell}}\right) +   N \langle D|Du|^2, D_{\bar{v}} \psi \rangle +  N \langle D_{\bar{v}} \psi, D|Du|^2 \rangle \\
&=& M \left(\langle Du, D_{\bar v} \psi \rangle+ \langle D_{\bar v} \psi, Du \rangle \right) \geq -CM.
\eean
Using (\ref{secondorder}), one can obtain the well-known identity
\be
-\sigma_k^{p \bar{q},r \bar{s}} D_j g_{\bar{q} p} D_{\bar{j}} g_{\bar{s} r} = -\sigma_k^{p \bar{p}, q \bar{q}} D_j g_{\bar{p} p} D_{\bar{j}} g_{\bar{q} q} + \sigma_k^{p \bar{p}, q \bar{q}} |D_j g_{\bar{p} q}|^2,
\ee
where $\sigma_k^{p \bar{p}, q \bar{q}} = {\p \over \p \lambda_p} {\p \over \p \lambda_q}  \sigma_k(\lambda)$. We assume that $\lambda_1 \gg 1$, otherwise the $C^2$ estimate is complete. The main inequality (\ref{mainineq}) becomes
\bea\label{mainineq2}
0 &\geq& {-C \over \lambda_1} \bigg\{ 1 + |DDu|^2 + |D \bar{D} u|^2 \bigg\} +  {1 \over P_m} \bigg\{ \sum_j \lambda_j^{m-1}  (  -\sigma_k^{p \bar{p}, q \bar{q}} D_j g_{\bar{p} p} D_{\bar{j}} g_{\bar{q} q} + \sigma_k^{p \bar{p}, q \bar{q}} |D_j g_{\bar{p} q}|^2 
\bigg\} \nonumber\\
&&+ {1 \over P_m} \bigg\{ (m-1)\sum_j \lambda_j^{m-2} \sigma_k^{p \bar{p}} |D_p g_{\bar{j} j}|^2 + \sigma_k^{p \bar{p}} \sum_{i \neq j} {\lambda_i^{m-1} - \lambda_j^{m-1} \over \lambda_i - \lambda_j} |D_p g_{\bar{j} i}|^2 \bigg\} \nonumber\\
&&- {|D P_m|^2_\sigma \over m P_m^2} + N (|DDu|^2_{\sigma \o} + |D \bar{D} u|^2_{\sigma \o}) + (M\e- C N -C) {\mathcal{F}} - CM.
\eea

\ 

The main objective is to show that the third order terms on the right hand side of (\ref{mainineq2}) are nonnegative. To deal with this issue, we need a lemma from \cite{GRW} (see also \cite{GuanLL, LRW}).

\begin{lemma}[\cite{GRW}] \label{GRW_ineq}
Suppose $1\leq \ell < k\leq n$, and let $\alpha = 1/(k-\ell)$. Let $W=(w_{\bar{q} p})$ be a Hermitian tensor in the $\Gamma_k$ cone. Then, for any $\theta>0$,
\bea
& \ &
-\sigma_k^{p \bar{p},q \bar{q}}(W) w_{\bar{p} p i} w_{\bar{q} q \bar{i}} + (1-\alpha + {\alpha \over \theta}) {|D_i \sigma_k(W)|^2 \over \sigma_k(W)} \nonumber\\
&\geq& \sigma_k(W) (\alpha+1-\alpha \theta) \bigg| { D_i \sigma_\ell(W) \over \sigma_\ell(W)} \bigg|^2 - {\sigma_k \over \sigma_\ell} (W) \sigma_\ell^{p \bar{p}, q \bar{q}}(W) w_{\bar{p} p i} w_{\bar{q} q \bar{i}} .
\eea
\end{lemma}
Here the $\Gamma_k$ cone is defined as following:
\bea\label{cone}
\Gamma_k = \ \{ \lambda \in \R^n \ | \ \sigma_m(\lambda)>0, \ m =1, \cdots, k\}.
\eea
We say a Hermitian matrix $W\in \Gamma_k$ if $\lambda(W) \in \Gamma_k$.

\bigskip

It follows from the above lemma that, by taking $\ell=1$, we have
\bea\label{ineq1}
-\sigma_k^{p \bar{p}, q \bar{q}} D_i g_{\bar{p} p} D_{\bar{i}} g_{\bar{q} q} + K  |D_i \sigma_k|^2 \geq 0,
\eea
for $K > (1-\alpha+ \frac{\alpha}{\theta}) \left(\inf \psi\right)^{-1}$ if $2\leq k\leq n$. 
\

\par
We shall denote
\bean
A_i = {\lambda_i^{m-1} \over P_m} \bigg\{ K |D_i \sigma_k|^2  -\sigma_k^{p \bar{p}, q \bar{q}} D_i g_{\bar{p} p} D_{\bar{i}} g_{\bar{q} q} \bigg\},
\eean
\bean
B_i = {1 \over P_m}  \bigg\{\sum_p \sigma_k^{p \bar{p}, i \bar{i}} \lambda_p^{m-1} |D_i g_{\bar{p} p}|^2  \bigg\},\ \ \  
C_i = {(m-1) \sigma_k^{i \bar{i}} \over P_m} \bigg\{ \sum_p \lambda_p^{m-2} |D_i g_{\bar{p} p}|^2 \bigg\},
\eean
\bean
D_i = {1 \over P_m} \bigg\{ \sum_{p \neq i} \sigma_k^{p \bar{p}} \ {\lambda_p^{m-1} - \lambda_i^{m-1} \over \lambda_p - \lambda_i} |D_i g_{\bar{p} p}|^2 \bigg\}, \ \ \ \ E_i  ={m \sigma_k^{i \bar{i}} \over P_m^2} \bigg| \sum_p \lambda_p^{m-1} D_i g_{\bar{p} p} \bigg|^2.
\eean
Define $T_{j \bar{p} q} = D_j \chi_{\bar{p} q} - D_q\chi_{\bar{p} j}$.
For any $0<\tau<1$, we can estimate
\bea
{1 \over P_m} \bigg\{ \sum_p \lambda_p^{m-1} \sigma_k^{j \bar{j}, i \bar{i}} |D_p g_{\bar{j} i}|^2 \bigg\} &\geq& {1 \over P_m} \bigg\{ \sum_p \lambda_p^{m-1} \sigma_k^{p \bar{p}, i \bar{i}} |D_i g_{\bar{p} p} + T_{p \bar{p} i}|^2 \bigg\}\nonumber\\
&\geq& {1 \over P_m} \bigg\{ \sum_p \lambda_p^{m-1} \sigma_k^{p \bar{p}, i \bar{i}} \{ (1-\tau) |D_i g_{\bar{p} p}|^2 - C_{\tau} |T_{p \bar{p} i}|^2 \} \bigg\}\nonumber\\
&=& (1-\tau) \sum_i B_i - {C_{\tau} \over P_m} \sum_p \lambda_p^{m-2} (\lambda_p \sigma_k^{p \bar{p}, i \bar{i}}) |T_{p \bar{p} i}|^2. \nonumber
\eea
Now, we use $\sigma_{l} (\lambda | i )$ and $\sigma_{l}(\lambda | ij )$ to denote the $l$-th elementary function of 
\bea
\nonumber(\lambda | i) = (\lambda_1, \cdots, \widehat{\lambda_{i}}, \cdots, \lambda_n)\in \R^{n-1} \textit{ and } \
(\lambda | ij) = (\lambda_1, \cdots, \widehat{\lambda_{i}}, \cdots, \widehat{\lambda_{j}}, \cdots, \lambda_n)\in \R^{n-2}
\eea
respectively. The following simple identities are used frequently,
\bea
\sigma_k^{i\bar i}=\sigma_{k-1}(\lambda| i),
\qquad
\sigma_k^{p\bar p,i\bar i}
=
\sigma_{k-2}(\lambda| p i).
\nonumber
\eea
Using the identity $\sigma_l (\lambda) = \sigma_{l} (\lambda | p) + \lambda_p \sigma_{l-1} (\lambda | p)$ for any $1\leq p\leq n$, we obtain
\bea\label{inequality0}
{1 \over P_m} \bigg\{ \sum_p \lambda_p^{m-1} \sigma_k^{j \bar{j}, i \bar{i}} |D_p g_{\bar{j} i}|^2 \bigg\}&\geq & (1-\tau) \sum_i B_i - {C_{\tau} \over P_m} \sum_p \lambda_p^{m-2} (\sigma_k^{i \bar{i}} - \sigma_{k-1}(\lambda | pi)) |T_{p \bar{p} i}|^2 \nonumber\\
&\geq& (1-\tau) \sum_i B_i - {C_{\tau} \over \lambda_1^2} {\mathcal{F}} \geq (1-\tau) \sum_i B_i - {\mathcal{F}}.
\eea
We used the notation $C_\tau$ for a constant depending on $\tau$. To get the last inequality above, we assumed that $\lambda_1^2\geq C_{\tau}$; otherwise, we already have the desired estimate $\lambda_1 \leq C$. Similarly, we may estimate
\bea
{1 \over P_m}\sigma_k^{j \bar{j}} \sum_{i \neq p} {\lambda_i^{m-1} - \lambda_p^{m-1} \over \lambda_i - \lambda_p} |D_j g_{\bar{p} i}|^2  &\geq& {1 \over P_m}\sigma_k^{p \bar{p}} \sum_{i; p \neq i} {\lambda_i^{m-1} - \lambda_p^{m-1} \over \lambda_i - \lambda_p} |D_i g_{\bar{p} p} + T_{p \bar{p} i}|^2\\
&\geq& {1 \over P_m}\sigma_k^{p \bar{p}} \sum_{i; p \neq i} {\lambda_i^{m-1} - \lambda_p^{m-1} \over \lambda_i - \lambda_p} \{ (1-\tau)|D_i g_{\bar{p} p}|^2 -C_{\tau} |T_{p \bar{p} i}|^2 \} \nonumber\\
&\geq& \sum_i (1-\tau) D_i - {C_{\tau} \over \lambda_1^2} {\mathcal{F}} \geq  \sum_i (1-\tau) D_i - \mathcal{F} \nonumber.
\eea
With the introduced notation in place, the main inequality becomes
\bea \label{inequality1}
0 &\geq& {-C(K) \over \lambda_1} \bigg\{ 1 + |DDu|^2 + |D \bar{D} u|^2 \bigg\} - \tau {|D P_m|^2_\sigma \over m P_m^2} \nonumber\\
&&+  \sum_i \bigg\{ A_i + (1-\tau)B_i +C_i +(1-\tau)D_i - (1-\tau)E_i \bigg\}  \nonumber\\
&&+ N (|DDu|^2_{\sigma \o} + |D \bar{D} u|^2_{\sigma \o}) + (M\e- C N -C) {\mathcal{F}} - CM.
\eea
Using the critical equation (\ref{testfunction1}), we have
\bea
\tau {|D P_m|^2_\sigma \over m P_m^2} &=& \tau m \bigg| N D |Du|^2 - M Du \bigg|^2_\sigma \leq 2 \tau m ( N^2 | D|Du|^2|^2_\sigma + M^2 |Du|^2_\sigma) \nonumber\\
&\leq& C \tau m N^2 (|DDu|^2_{\sigma \o} + |D \bar{D} u|^2_{\sigma \o}) + C \tau m M^2 \mathcal{F}.
\eea
We thus have
\bea \label{inequality_before_3rd_order_est}
0 &\geq& {-C(K) \over \lambda_1} \bigg\{ 1 + |DDu|^2 + |D \bar{D} u|^2 \bigg\}  + (N-C\tau m N^2) (|DDu|^2_{\sigma \o} + |D \bar{D} u|^2_{\sigma \o}) \nonumber\\
&&+\sum_i \bigg\{ A_i + (1-\tau)B_i +C_i +(1-\tau)D_i - (1-\tau)E_i \bigg\} \nonumber\\
&& + (M\e - C \tau m M^2- C N -C) {\mathcal{F}} - CM. 
\eea

\subsection{Estimating the Third Order Terms}
In this subsection, we will adapt the argument in \cite{LRW} to estimate the third order terms.
\begin{lemma} \label{Cauchy-Schwartz} For sufficiently large $m$, the following estimates hold:
\bea \label{genesis2}
P_m^2 (B_1 + C_1 +D_1 - E_1) \geq P_m \lambda_1^{m-2} \sum_{p \neq 1} \sigma_k^{p \bar{p}} |D_1 g_{\bar{p} p}|^2 - \lambda_1^m \sigma_k^{1 \bar{1}}  \lambda_1^{m-2} |D_1 g_{\bar{1} 1}|^2,
\eea
and for any fixed $i \neq 1$,
\be\label{genesisi}
P_m^2 (B_i + C_i + D_i -E_i) \geq 0.
\ee
\end{lemma}
{\it Proof.} Fix $i \in \{ 1,2, \dots, n \}$. First, we compute
\bea
P_m (B_i + D_i) &=& \sum_{p \neq i} \sigma_k^{p \bar{p}, i \bar{i}} \lambda_p^{m-1} |D_i g_{\bar{p} p}|^2  +  \sum_{p \neq i} \sigma_k^{p \bar{p}} \ {\lambda_p^{m-1} - \lambda_i^{m-1} \over \lambda_p - \lambda_i} |D_i g_{\bar{p} p}|^2 \nonumber\\
&=& \sum_{p \neq i} \lambda_p^{m-2} \bigg\{  (\lambda_p \sigma_k^{p \bar{p}, i \bar{i}} +  \sigma_k^{p \bar{p}}) |D_i g_{\bar{p} p}|^2\bigg\} +\bigg\{ \sum_{p \neq i} \sigma_k^{p \bar{p}} \sum_{q=0}^{m-3} \lambda_p{}^q \lambda_i^{m-2-q} |D_i g_{\bar{p} p}|^2 \bigg\}. \nonumber
\eea
Note that, 
\bea
\lambda_p\sigma_k^{p\bar p,i\bar i}+\sigma_k^{p\bar p}\geq \sigma_k^{i\bar i}.
\nonumber
\eea
To see this, we write
\bea
\lambda_p\sigma_k^{p\bar p,i\bar i}+\sigma_k^{p\bar p}&=&\lambda_p \sigma_{k-2}(\lambda|pi)+
\sigma_{k-1}(\lambda|p)
\nonumber\\
&=&\sigma_{k-1}(\lambda|i)-\sigma_{k-1}(\lambda|ip)+\sigma_{k-1}(\lambda|p)
\nonumber\\
&=&\sigma_{k-1}(\lambda|i)+\lambda_i\sigma_{k-2}(\lambda|ip)
\geq \sigma_{k-1}(\lambda|i)=\sigma_k^{i\bar i},
\nonumber
\eea
where we used the standard identity $\sigma_l (\lambda) = \sigma_{l} (\lambda | p) + \lambda_p \sigma_{l-1} (\lambda | p)$ twice, to get the second and third equalities.
Therefore
\be
P_m (B_i + D_i) \geq \sigma_k^{i \bar{i}} \bigg\{ \sum_{p \neq i} \lambda_p^{m-2} |D_i g_{\bar{p} p}|^2\bigg\} +\bigg\{ \sum_{p \neq i} \sigma_k^{p \bar{p}} \sum_{q=0}^{m-3} \lambda_p{}^q \lambda_i^{m-2-q} |D_i g_{\bar{p} p}|^2 \bigg\}.
\ee
It follows that 
\bea
P_m (B_i + C_i + D_i) &\geq& m \sigma_k^{i \bar{i}}  \sum_{p \neq i} \lambda_p^{m-2} |D_i g_{\bar{p} p}|^2 + (m-1)\sigma_k^{i \bar{i}}  \lambda_i^{m-2} |D_i g_{\bar{i} i}|^2  \nonumber\\
&&+  \sum_{p \neq i} \sigma_k^{p \bar{p}} \sum_{q=0}^{m-3} \lambda_p{}^q \lambda_i^{m-2-q} |D_i g_{\bar{p} p}|^2.
\eea
Expanding out the definition of $E_i$
\be
P_m^2 E_i = m \sigma_k^{i \bar{i}} \sum_{p \neq i} \lambda_p^{2m -2} |D_i g_{\bar{p} p}|^2 + m \sigma_k^{i \bar{i}} \lambda_i^{2m-2} |D_i g_{\bar{i} i}|^2 + m \sigma_k^{i \bar{i}} \sum_p \sum_{q \neq p} \lambda_p^{m-1} \lambda_q^{m-1} D_i g_{\bar{p} p} D_{\bar{i}} g_{\bar{q} q}.
\ee
Therefore
\bea
&& P_m^2 (B_i + C_i + D_i - E_i)\\ &\geq&\bigg\{ m \sigma_k^{i \bar{i}}  \sum_{p \neq i} (P_m - \lambda_p^m) \lambda_p^{m-2} |D_i g_{\bar{p} p}|^2 - m \sigma_k^{i \bar{i}} \sum_{p \neq i} \sum_{q \neq p,i} \lambda_p^{m-1} \lambda_q^{m-1} D_i g_{\bar{p} p} D_{\bar{i}} g_{\bar{q} q} \bigg\} \nonumber\\
&&+ P_m \sum_{p \neq i} \sigma_k^{p \bar{p}} \sum_{q=0}^{m-3} \lambda_p{}^q \lambda_i^{m-2-q} |D_i g_{\bar{p} p}|^2 - 2 m \sigma_k^{i \bar{i}} \Re \sum_{q \neq i} \lambda_i^{m-1} \lambda_q^{m-1} D_i g_{\bar{i} i} D_{\bar{i}} g_{\bar{q} q} \nonumber\\
&&+ \{(m-1)P_m - m \lambda_i^m\} \sigma_k^{i \bar{i}}  \lambda_i^{m-2} |D_i g_{\bar{i} i}|^2.\nonumber
\eea
We shall estimate the expression in brackets. First, 
\bean
m \sigma_k^{i \bar{i}}  \sum_{p \neq i} (P_m - \lambda_p^m) \lambda_p^{m-2} |D_i g_{\bar{p} p}|^2 = m \sigma_k^{i \bar{i}}  \sum_{p \neq i} \sum_{q \neq p,i} \lambda_q^m \lambda_p^{m-2} |D_i g_{\bar{p} p}|^2 +  m \sigma_k^{i \bar{i}}  \sum_{p \neq i} \lambda_i^m \lambda_p^{m-2} |D_i g_{\bar{p} p}|^2.
\eean
Next, we can estimate
\bea
&&-m \sigma_k^{i \bar{i}} \sum_{p \neq i} \sum_{q \neq p,i} \lambda_p^{m-1} \lambda_q^{m-1} D_i g_{\bar{p} p} D_{\bar{i}} g_{\bar{q} q} \\ &\geq& -m \sigma_k^{i \bar{i}} \sum_{p \neq i} \sum_{q \neq p,i} {1 \over 2} \{\lambda_p^{m-2} \lambda_q^{m}  |D_i g_{\bar{p} p}|^2 + \lambda_p^{m} \lambda_q^{m-2}  |D_i g_{\bar{q} q}|^2\} = -m \sigma_k^{i \bar{i}} \sum_{p \neq i} \sum_{q \neq p,i} \lambda_p^{m-2} \lambda_q^{m}  |D_i g_{\bar{p} p}|^2.\nonumber
\eea
We arrive at
\bea \label{genesis1}
&&P_m^2 (B_i + C_i + D_i - E_i)\\ &\geq& m \sigma_k^{i \bar{i}}  \sum_{p \neq i} \lambda_i^m \lambda_p^{m-2} |D_i g_{\bar{p} p}|^2 + P_m \sum_{p \neq i} \sigma_k^{p \bar{p}} \sum_{q=0}^{m-3} \lambda_p{}^q \lambda_i^{m-2-q} |D_i g_{\bar{p} p}|^2 \nonumber\\
&& - 2 m \sigma_k^{i \bar{i}}  \Re \bigg\{ \lambda_i^{m-1} D_i g_{\bar{i} i} \sum_{q \neq i} \lambda_q^{m-1}  D_{\bar{i}} g_{\bar{q} q} \bigg\} + \{(m-1)P_m - m \lambda_i^m\} \sigma_k^{i \bar{i}}  \lambda_i^{m-2} |D_i g_{\bar{i} i}|^2.\nonumber
\eea
The next step is to extract good terms from the second summation on the first line. We fix a $p\neq i$. \newline
\newline
Case 1: $\lambda_i \geq \lambda_p$. Then $\sigma_k^{p \bar{p}} \geq \sigma_k^{i \bar{i}}$. Hence
\be
P_m \sigma_k^{p \bar{p}} \sum_{q=1}^{m-3} \lambda_p{}^q \lambda_i^{m-2-q} \geq \lambda_i^m \sigma_k^{i \bar{i}} \sum_{q=1}^{m-3} \lambda_p{}^q \lambda_p^{m-2-q} = (m-3) \sigma_k^{i \bar{i}} \lambda_i^m \lambda_p^{m-2}.
\ee
Case 2: $\lambda_i \leq \lambda_p$. Then $\lambda_p \sigma_k^{p \bar{p}} = \lambda_i \sigma_k^{i \bar{i}} + (\sigma_k(\lambda|i) - \sigma_k(\lambda|p)) \geq \lambda_i \sigma_k^{i \bar{i}}$,
and we obtain
\be
P_m \sigma_k^{p \bar{p}} \sum_{q=1}^{m-3} \lambda_p{}^q \lambda_i^{m-2-q} \geq \lambda_p^m \sigma_k^{i \bar{i}} \sum_{q=1}^{m-3} \lambda_p{}^{q-1} \lambda_i^{m-1-q} \geq (m-3) \sigma_k^{i \bar{i}} \lambda_i^m \lambda_p^{m-2}.
\ee
Combining both cases, we have
\bea
P_m \sigma_k^{p \bar{p}} \sum_{q=0}^{m-3} \lambda_p{}^q \lambda_i^{m-2-q} |D_i g_{\bar{p} p}|^2 &=& P_m \sigma_k^{p \bar{p}} \sum_{q=1}^{m-3} \lambda_p{}^q \lambda_i^{m-2-q} |D_i g_{\bar{p} p}|^2 + P_m \sigma_k^{p \bar{p}} \lambda_i^{m-2} |D_i g_{\bar{p} p}|^2 \nonumber\\
&\geq&  (m-3) \sigma_k^{i \bar{i}} \lambda_i^m \lambda_p^{m-2}|D_i g_{\bar{p} p}|^2 + P_m \sigma_k^{p \bar{p}} \lambda_i^{m-2} |D_i g_{\bar{p} p}|^2. \nonumber
\eea
Substituting this estimate into inequality (\ref{genesis1}), we obtain
\bea 
&& P_m^2 (B_i + C_i + D_i - E_i)\\
 &\geq& (2m-3) \sigma_k^{i \bar{i}}  \sum_{p \neq i} \lambda_i^m \lambda_p^{m-2} |D_i g_{\bar{p} p}|^2 - 2 m \sigma_k^{i \bar{i}}  \Re \bigg\{ \lambda_i^{m-1} D_i g_{\bar{i} i} \sum_{p \neq i} \lambda_p^{m-1}  D_{\bar{i}} g_{\bar{p} p} \bigg\} \nonumber\\
&& + P_m \lambda_i^{m-2} \sum_{p \neq i} \sigma_k^{p \bar{p}} |D_i g_{\bar{p} p}|^2 + \{(m-1)P_m - m \lambda_i^m\} \sigma_k^{i \bar{i}}  \lambda_i^{m-2} |D_i g_{\bar{i} i}|^2.\nonumber
\eea
Choose $m \gg 1$ such that
\be
m^2 \leq (2m-3)(m-2).
\ee
We can therefore estimate
\bea 
&\ & 2 m \sigma_k^{i \bar{i}} \Re \bigg\{ \lambda_i^{m-1} D_i g_{\bar{i} i} \sum_{p \neq i} \lambda_p^{m-1}  D_{\bar{i}} g_{\bar{p} p} \bigg\} \nonumber\\
&\leq&  2  \sigma_k^{i \bar{i}} \sum_{p \neq i} \{ (2m-3)^{1/2} \lambda_i^{m/2} \lambda_p^{m-2 \over 2} |D_i g_{\bar{p} p}| \} \{ (m-2)^{1/2}  \lambda_i^{m-2 \over 2} \lambda_p^{m/2} |D_{\bar{i}} g_{\bar{i} i}|\} \nonumber\\
&\leq& (2m-3) \sigma_k^{i \bar{i}} \sum_{p \neq i}  \lambda_i^{m} \lambda_p^{m-2} |D_i g_{\bar{p} p}|^2 + (m-2) \sigma_k^{i \bar{i}} \sum_{p \neq i} \lambda_i^{m-2} \lambda_p^{m} |D_{\bar{i}} g_{\bar{i} i}|^2.
\eea
We finally arrive at
\bea
P_m^2 (B_i + C_i +D_i - E_i) &\geq& P_m \lambda_i^{m-2} \sum_{p \neq i} \sigma_k^{p \bar{p}} |D_i g_{\bar{p} p}|^2 + \{(m-1)P_m - m \lambda_i^m\} \sigma_k^{i \bar{i}}  \lambda_i^{m-2} |D_i g_{\bar{i} i}|^2 \nonumber\\
&&-(m-2) \sigma_k^{i \bar{i}} \sum_{p \neq i} \lambda_i^{m-2} \lambda_p^{m} |D_{\bar{i}} g_{\bar{i} i}|^2.
\eea
If we let $i=1$, we obtain inequality (\ref{genesis2}). For any fixed $i \neq 1$, this inequality yields
\bea
P_m^2 (B_i + C_i +D_i - E_i) &\geq& P_m \lambda_i^{m-2} \sum_{p \neq i} \sigma_k^{p \bar{p}} |D_i g_{\bar{p} p}|^2 + \{(m-1)\lambda_1^m - \lambda_i^m\} \sigma_k^{i \bar{i}}  \lambda_i^{m-2} |D_i g_{\bar{i} i}|^2 \nonumber\\
&&+ (m-1) \sum_{p \neq 1,i} \lambda_p^m \sigma_k^{i \bar{i}}  \lambda_i^{m-2} |D_i g_{\bar{i} i}|^2 -(m-2) \sigma_k^{i \bar{i}} \sum_{p \neq i} \lambda_i^{m-2} \lambda_p^{m} |D_{\bar{i}} g_{\bar{i} i}|^2 \nonumber\\
&\geq& P_m \lambda_i^{m-2} \sum_{p \neq i} \sigma_k^{p \bar{p}} |D_i g_{\bar{p} p}|^2 \geq 0. \nonumber
\eea
This completes the proof of Lemma \ref{Cauchy-Schwartz}. Q.E.D.

\bigskip

We observed in (\ref{ineq1}) that $A_i\geq 0$. Lemma \ref{Cauchy-Schwartz} implies that for any $i\neq 1$, 
\bea
A_i + B_i + C_i + D_i -E_i \geq 0. \nonumber
\eea
Thus we have shown that for $i\neq 1$, the third order terms in the main inequality (\ref{inequality_before_3rd_order_est}) are indeed nonnegative. The only remaining case is when $i=1$. 
By adapting once again the techniques from \cite{GRW}, we obtain the following lemma.

\begin{lemma}\label{key_lemma}
Let $1<k\leq n$. Suppose there exists $0<\delta \leq 1$ such that $\lambda_\mu \geq \delta \lambda_1$ for some $\mu \in \{ 1,2, \dots, k-1 \}$. There exists a small $\delta'>0$ such that if $\lambda_{\mu+1} \leq \delta' \lambda_1$, then
\bea
A_1 + B_1 + C_1 + D_1 -E_1 \geq 0. \nonumber
\eea
\end{lemma}
{\it Proof.} By Lemma \ref{Cauchy-Schwartz}, we have
\bea \label{genesis3}
&&P_m^2 (A_1+ B_1 + C_1 +D_1 - E_1) \\
&\geq&
P_m^2 A_1 + P_m \lambda_1^{m-2} \sum_{p \neq 1} \sigma_k^{p \bar{p}} |D_1 g_{\bar{p} p}|^2 - \lambda_1^m \sigma_k^{1 \bar{1}}  \lambda_1^{m-2} |D_1 g_{\bar{1} 1}|^2. \nonumber
\eea
The key insight in \cite{GRW}, used also in \cite{LRW}, is to extract a good term involving $|D_1 g_{\bar{1} 1}|^2$ from $A_1$. By the inequality in Lemma \ref{GRW_ineq} (with $\theta = {1\over 2}$), we have for $\mu < k$
\bea \label{A_first_estimate}
P_m^2 A_1 &\geq& {P_m \lambda_1^{m-1} \sigma_k \over \sigma_\mu^2 } \bigg\{  (1+{\alpha \over 2}) \bigg| \sum_p \sigma_\mu^{p \bar{p}} D_1 g_{\bar{p} p} \bigg|^2 - \sigma_\mu \sigma_\mu^{p \bar{p}, q \bar{q}} D_1 g_{\bar{p} p} D_{\bar{1}} g_{\bar{q} q} \bigg\} \nonumber\\
&=& {P_m \lambda_1^{m-1} \sigma_k \over \sigma_\mu^2 } \bigg\{ \sum_{p} (1+{  \alpha \over 2} )|\sigma_\mu^{p \bar{p}} D_1 g_{\bar{p} p}|^2 + \sum_{p \neq q} {\alpha \over 2}\sigma_\mu^{p \bar{p}} D_1 g_{\bar{p} p} \sigma_\mu^{q \bar{q}} D_{\bar{1}} g_{\bar{q} q}  \nonumber\\
&&+ \sum_{p \neq q} (\sigma_\mu^{p \bar{p}} \sigma_\mu^{q \bar{q}} - \sigma_\mu \sigma_\mu^{p \bar{p}, q \bar{q}}) D_1 g_{\bar{p} p} D_{\bar{1}} g_{\bar{q} q} \bigg\} \nonumber\\
&\geq&  {P_m \lambda_1^{m-1} \sigma_k \over \sigma_\mu^2 } \bigg\{ \sum_p |\sigma_\mu^{p \bar{p}} D_1 g_{\bar{p} p}|^2 - \sum_{p \neq q} | F^{pq} D_1 g_{\bar{p} p} D_{\bar{1}} g_{\bar{q} q} |\bigg\},
\eea
where we defined $F^{pq} = \sigma_\mu^{p \bar{p}} \sigma_\mu^{q \bar{q}} - \sigma_\mu \sigma_\mu^{p \bar{p}, q \bar{q}}$. 
Notice if $\mu = 1$, then $F^{pq} = 1$. If $\mu \geq 2$, then the Newton-MacLaurin inequality implies
\be \label{pq_identity}
F^{pq} = \sigma_{\mu-1}^2(\lambda|pq) - \sigma_\mu(\lambda|pq) \sigma_{\mu-2}(\lambda|pq) \geq 0.
\ee
 We split the sum involving $F^{pq}$ in the following way:
\bea \label{split_sum}
\sum_{p \neq q} |F^{pq} D_1 g_{\bar{p} p} D_{\bar{1}} g_{\bar{q} q}| &=& \sum_{p \neq q; p,q \leq \mu} F^{pq} |D_1 g_{\bar{p} p}| |D_{\bar{1}} g_{\bar{q} q}| + \sum_{(p,q) \in J} F^{pq} |D_1 g_{\bar{p} p}| |D_{\bar{1}} g_{\bar{q} q}|
\eea
where $J$ is the set of indices where at least one of $p \neq q$ is strictly greater than $\mu$. The summation of terms in $J$ can be estimated by
\bea 
- \sum_{(p,q) \in J} F^{pq} |D_1 g_{\bar{p} p}| |D_{\bar{1}} g_{\bar{q} q}| &\geq&- \sum_{(p,q) \in J} \sigma_\mu^{p \bar{p}} \sigma_\mu^{q \bar{q}} |D_1 g_{\bar{p} p}| |D_{\bar{1}} g_{\bar{q} q}| \nonumber\\
&\geq& - \epsilon \sum_{p \leq \mu} |\sigma_\mu^{p \bar{p}} D_1 g_{\bar{p} p}|^2 -C \sum_{p > \mu} |\sigma_\mu^{p \bar{p}} D_1 g_{\bar{p} p}|^2.
\eea
\par
If $\mu =1$, the first term on the right hand side of (\ref{split_sum}) vanishes and this estimate applies to all terms on the right hand side of (\ref{split_sum}). 

If $\mu \geq 2$, we have for $p,q \leq \mu$, 
\be \label{mu-1}
\sigma_{\mu-1}(\lambda|pq) \leq C {\lambda_1 \cdots \lambda_{\mu+1} \over \lambda_p \lambda_q} \leq C {\sigma_\mu^{p \bar{p}} \lambda_{\mu+1} \over \lambda_q}.
\ee
Using (\ref{pq_identity}) and (\ref{mu-1}), for $\delta'$ small enough we can control
\bea
&&- \sum_{p \neq q; p,q \leq \mu} F^{pq} |D_1 g_{\bar{p} p}| |D_{\bar{1}} g_{\bar{q} q}| \geq - \sum_{p \neq q; p,q \leq \mu} \sigma_{\mu-1}^2(\lambda|pq) |D_1 g_{\bar{p} p}| |D_{\bar{1}} g_{\bar{q} q}| \nonumber\\
&\geq& - C\lambda_{\mu+1}^2\sum_{p \neq q; p,q \leq \mu}  {\sigma_\mu^{p \bar{p}} \over \lambda_p} |D_1 g_{\bar{p} p}|{\sigma_\mu^{q \bar{q}} \over \lambda_q}  |D_{\bar{1}} g_{\bar{q} q}| \geq - C \sum_{p \leq \mu} {\lambda_{\mu+1}^2 \over \lambda_p^2} |\sigma_\mu^{p \bar{p}} D_1 g_{\bar{p} p}|^2 \nonumber\\
&\geq& - C \sum_{p \leq \mu} {\delta'^2 \over \delta^2} |\sigma_\mu^{p \bar{p}} D_1 g_{\bar{p} p}|^2 \geq - \epsilon \sum_{p \leq \mu} |\sigma_\mu^{p \bar{p}} D_1 g_{\bar{p} p}|^2.
\eea
Combining all cases, we have
\be
- \sum_{p \neq q} | F^{pq} D_1 g_{\bar{p} p} D_{\bar{1}} g_{\bar{q} q} | \geq - 2 \epsilon \sum_{p \leq \mu} |\sigma_\mu^{p \bar{p}} D_1 g_{\bar{p} p}|^2 -C \sum_{p > \mu} |\sigma_\mu^{p \bar{p}} D_1 g_{\bar{p} p}|^2.
\ee
Using this inequality in (\ref{A_first_estimate}) yields 
\bea
P_m^2 A_1 &\geq& {P_m \lambda_1^{m-1} \sigma_k \over \sigma_\mu^2 } \bigg\{ (1-2 \epsilon) \sum_{p \leq \mu} |\sigma_\mu^{p \bar{p}} D_1 g_{\bar{p} p}|^2 -C \sum_{p > \mu} |\sigma_\mu^{p \bar{p}} D_1 g_{\bar{p} p}|^2  |\bigg\} \nonumber\\
&\geq& (1-2 \epsilon) {P_m \lambda_1^{m-1} \sigma_k \over \sigma_\mu^2 } |\sigma_\mu^{1 \bar{1}} D_1 g_{\bar{1} 1}|^2 - C {P_m \lambda_1^{m-1} \sigma_k \over \sigma_\mu^2 } \sum_{p > \mu} |\sigma_\mu^{p \bar{p}} D_1 g_{\bar{p} p}|^2.
\eea
We estimate
\bea
&&(1 -2 \epsilon){P_m \lambda_1^{m-1} \sigma_k \over \sigma_\mu^2 } |\sigma_\mu^{1 \bar{1}} D_1 g_{\bar{1} 1}|^2 = (1 - 2\epsilon){P_m \lambda_1^{m-2} \sigma_k \over \lambda_1 } \bigg( {\lambda_1 \sigma_\mu^{1 \bar{1}} \over \sigma_\mu} \bigg)^2 |D_1 g_{\bar{1} 1}|^2 \nonumber\\ 
&\geq& (1 - 2\epsilon )P_m \lambda_1^{m-2} {\sigma_k \over \lambda_1 } \bigg( 1 - C{\lambda_{\mu+1} \over \lambda_1} \bigg)^2 |D_1 g_{\bar{1} 1}|^2  \geq (1 - 2 \epsilon)(1 - C \delta')^2 P_m \lambda_1^{m-2} \sigma_k^{1 \bar{1}} |D_1 g_{\bar{1} 1}|^2 \nonumber\\
&\geq& (1 - 2 \epsilon)(1 - C \delta')^2 (1+\delta^m) \lambda_1^{2m-2} \sigma_k^{1 \bar{1}} |D_1 g_{\bar{1} 1}|^2.
\eea
For $\delta'$ and $\epsilon$ small enough, we obtain
\be
P_m^2 A_1 \geq  \lambda_1^m \sigma_k^{1 \bar{1}} \lambda_1^{m-2} |D_1 g_{\bar{1} 1}|^2 - C {P_m \lambda_1^{m-1} \sigma_k \over \sigma_\mu^2 } \sum_{p > \mu} |\sigma_\mu^{p \bar{p}} D_1 g_{\bar{p} p}|^2.
\ee
We see that the $|D_1 g_{\bar{1} 1}|^2$ term cancels from inequality (\ref{genesis3}) and we are left with
\be
P_m^2 (A_1 + B_1 + C_1 +D_1 - E_1) \geq P_m \lambda_1^{m-2} \sum_{p > \mu} \bigg\{ \sigma_k^{p \bar{p}}- C {\lambda_1 \sigma_k (\sigma_\mu^{p \bar{p}})^2 \over \sigma_\mu^2} \bigg\} |D_1 g_{\bar{p} p}|^2.
\ee
For $\delta'$ small enough, the above expression is nonnegative. Indeed, for any $p >\mu$, we have
\be
(\lambda_1 \sigma_\mu^{p \bar{p}})^2 \leq {1 \over \delta^2} (\lambda_\mu \sigma_\mu^{p \bar{p}})^2 \leq C {(\sigma_\mu)^2 \over \delta^2},
\ee
Therefore
\be
C {\lambda_1 \sigma_k (\sigma_\mu^{p \bar{p}})^2 \over \sigma_\mu^2} \leq {C \over \delta^2} {\sigma_k \over \lambda_1}.
\ee
On the other hand,  we notice that, if $p>k$, then $\sigma_k^{p \bar{p}} \geq \lambda_1 \cdots \lambda_{k-1} \geq c_n {\sigma_k \over \lambda_k} \geq {c_n \over \delta'} {\sigma_k \over \lambda_1}$.
If $\mu<p \leq k$, then $\sigma_k^{p \bar{p}} \geq {\lambda_1 \cdots \lambda_{k} \over \lambda_p} \geq c_n {\sigma_k \over \lambda_p} \geq {c_n \over \delta'} {\sigma_k \over \lambda_1}$. 
%
It follows that for $\delta'$ small enough we have
\be \label{leftover_estimate}
\sigma_k^{p \bar{p}} \geq C {\lambda_1 \sigma_k (\sigma_\mu^{p \bar{p}})^2 \over \sigma_\mu^2}.
\ee
This completes the proof of Lemma \ref{key_lemma}. Q.E.D.

\subsection{Completing the Proof}
With Lemma \ref{Cauchy-Schwartz} and Lemma \ref{key_lemma} at our disposal, we claim that we may assume in inequality (\ref{inequality_before_3rd_order_est}) that
\be \label{claim1}
A_i + B_i + C_i + D_i - E_i \geq 0, \ \ \ \forall i=1, \cdots, n.
\ee
Indeed, first set $\delta_1 =1 $. If $\lambda_2 \leq \delta_2 \lambda_1$ for $\delta_2>0$ small enough, then by Lemma \ref{key_lemma} we see that (\ref{claim1}) holds. Otherwise, $\lambda_2 \geq \delta_2 \lambda_1$. If $\lambda_3 \leq \delta_3 \lambda_1$ for $\delta_3>0$ small enough, then by Lemma \ref{key_lemma} we see that (\ref{claim1}) holds. Otherwise, $\lambda_3 \geq \delta_3 \lambda_1$. Proceeding iteratively, we may arrive at $\lambda_k \geq \delta_k \lambda_1$. But in this case, the $C^2$ estimate follows directly from the equation as
\be\label{directly}
C \geq \sigma_k \geq \lambda_1 \cdots \lambda_k \geq (\delta_k)^{k-1} \lambda_1.
\ee
Therefore we may assume (\ref{claim1}), and inequality (\ref{inequality_before_3rd_order_est}) becomes
\bea 
0 &\geq& {-C(K) \over \lambda_1} \bigg\{ 1 + |DDu|^2 + |D \bar{D} u|^2 \bigg\}  + (N-C\tau m N^2) (|DDu|^2_{\sigma \o} + |D \bar{D} u|^2_{\sigma \o}) \nonumber\\
&&+ (M\e - C \tau m M^2- C N -C) {\mathcal{F}} - CM.
\eea
Since for fixed $i$, $\sigma_k^{i \bar{i}} \geq \sigma_k^{1 \bar{1}} \geq {k \over n} {\sigma_k \over \lambda_1} \geq {1 \over C \lambda_1}$,
we can estimate
\be
|DDu|^2_{\sigma \o} + |D \bar{D} u|^2_{\sigma \o} \geq {1 \over C \lambda_1}(|DDu|^2 + |D \bar{D} u|^2) \geq {1 \over C \lambda_1}|DDu|^2 + {\lambda_1 \over C} .
\ee
This leads to
\bean
0 &\geq& \bigg\{ {N \over C} - C \tau m N^2 -C(K) \bigg\} \lambda_1 +{1 \over \lambda_1} \bigg\{ {N \over C} - C \tau m N^2 -C(K) \bigg\} \bigg\{ 1 + |DDu|^2 \bigg\} \nonumber\\
&& + (M\e - C \tau m M^2- C N -C) {\mathcal{F}}- CM. \nonumber
\eean
By choosing $\tau$ small, for example, $\tau = \frac{1}{NM}$, we have
\bean
0 &\geq& \bigg\{ {N \over C} - {Cm \over M}N -C(K) \bigg\} \lambda_1+{1 \over \lambda_1} \bigg\{ {N \over C} - {C m \over M}N -C(K) \bigg\} \bigg\{ 1 + |DDu|^2 \bigg\} \nonumber\\
&&+ (M\e - {C m\over N} M- C N -C) {\mathcal{F}}  - CM. \nonumber
\eean
Taking $N$ and $M$ large enough, we can make the coefficients of the first three terms to be positive. For example, if we let $M=N^2$ for $N$ large, then ${N \over C} - {Cm \over M}N -C(K) = {N\over C} - {Cm\over N} - C(K)>0$ and $M\e - {C m\over N} M- C N -C= N^2\e-CmN-CN-C>0$. Thus, an upper bound of $\lambda_1$ follows. Q.E.D.

\begin{remark}\label{remarkcone}
In the above estimate, we assume that $\lambda=(\lambda_1, \cdots, \lambda_n)\in \Gamma_n$. Indeed, our estimate still works with $\lambda\in \Gamma_{k+1}$. It was observed in \cite{LRW} (Lemma 7) that if $\lambda \in \Gamma_{k+1}$, then $\lambda_1\geq \cdots \geq \lambda_n > -K_0$ for some positive constant $K_0$. Thus, we can replace $\lambda$ by $\tilde \lambda= \lambda+ K_0 I$ in our test function $G$ in (\ref{testfunction}).
\end{remark}

\bigskip

\noindent {\bf Acknowledgements:} The authors would like to thank Pengfei Guan for stimulating conversations and for sharing his unpublished notes \cite{Guan}. The authors are also very grateful to the referees for an exceptionally careful reading and for many helpful suggestions.

\bigskip
Department of Mathematics, Columbia University, New York, NY 10027, USA

\smallskip

phong@math.columbia.edu,  picard@math.columbia.edu, xzhang@math.columbia.edu

\bigskip
Department of Mathematics, University of California, Irvine, CA 92697, USA

\smallskip
xiangwen@math.uci.edu


\begin{thebibliography}{99}

{\small

\bibitem{AV} Alekser, S. and Verbitsky, M., {\it Quaternionic Monge-Amp\`ere equations and Calabi problem for HKT-manifolds}, Israel J. Math. 176 (2010), 109-138.

\bibitem{Ball} Ball, J., {\it Differentiability properties of symmetric and isotropic functions}, Duke Maht. J. 51 (1984), no. 3, 699-728.

\bibitem{Blocki1} Blocki, Z., {\it Weak solutions to the complex Hessian equation}, Ann. Inst. Fourier (Grenoble) 55, 5 (2005), 1735-1756.

\bibitem{DK} Dinew, S. and Kolodziej, S., {\it A priori estimates for the complex Hessian equations}, Analysis and PDE, Vol. 7, No. 1, (2014), 227-244.

\bibitem{FY} Fu, J.X. and Yau, S.T. {\it The theory of superstring with flux on non-K\"ahler manifolds and the complex Monge-Ampere equation}, J. Differential Geom, Vol 78, No. 3 (2008), 369-428.

\bibitem{FY1} Fu, J.X. and Yau, S.T. {\it A Monge-Amp\`ere type equation motivated by string theory}, Comm. in Geometry and Analysis, Vol 15, Number 1, (2007), 29-76.

\bibitem{Guanbo} Guan, B., {\it Second order estimates and regularity for fully nonlinear elliptic equations on Riemannian manifolds}, Duke Math. J. 163 (2014), 1491-1524.

\bibitem{Guan} Guan, P. and Ma, X.N., unpublished notes.

\bibitem{GuanLL} Guan, P., Li, J.F. and Li Y.Y., {\it Hypersurfaces of prescribed curvature measures}, 
Duke Math. J. Vol. 161, No. 10 (2012), 1927-1942.

\bibitem{GRW} Guan, P., Ren, C. and Wang, Z., {\it Global $C^2$ estimates for convex solutions of curvature equations}, Comm. Pure Appl. Math 68 (2015), 1287-1325.

\bibitem{Hou} Hou, Z., {\it Complex Hessian equation on K\"ahler manifold}, Int. Math. Res. Not. 16 (2009), 3098-3111.

\bibitem{HMW} Hou, Z., Ma, X.N. and Wu, D. {\it A second order estimate for complex Hessian equations on a compact K\"ahler manifold}, Math. Res. Lett. 17 (2010), 547-561.

\bibitem{KN} Kolodziej, S. and Nguyen, N-C., {\it Weak solutions of complex Hessian equations on compact Hermitian manifolds}, arXiv:1507.06755, to appear in Compos. Math.

\bibitem{LRW} Li, M., Ren, C. and Wang, Z., {\em An interior estimate for convex solutions and a rigidity theorem}, 
J. Funct. Anal., Vol. 270, Issue 7, (2016), 2691-2714.


\bibitem{Li} Li, S.Y., {\it On the Dirichlet problems for symmetric function equations of the eigenvalues of the complex Hessian}, Asian J. Math. 8 (2004), 87-106.

\bibitem{LN} Lu, H-C. and Nguyen, V-D., {\it Degenerate complex Hessian equations on compact K\"ahler manifolds}, preprint, arXiv: 1402.5147. to appear in Indiana Univ. Math. J.

\bibitem{PPZ} Phong, D.H., Picard, S. and Zhang, X.W., {\it On estimates for the Fu-Yau generalization of a Strominger system}, 	arXiv:1507.08193.

\bibitem{SW} Song, J. and Weinkove, B., {\it On the convergence and singularities of the J-flow with applications to the Mabuchi energy}, Comm. Pure. Appl. Math. 61 (2008), 210-229.

\bibitem{SX} Spruck, J. and Xiao, L., {\it A note on starshaped compact hypersurfaces with a prescribed scalar curvature in space forms}, arXiv:1505.01578.

\bibitem{Sun} Sun, W., {\it On uniform estimate of complex elliptic equations on closed Hermitian manifolds}, arXiv:1412.5001.

\bibitem{Gabor} Sz\'ekelyhidi, G., {\it Fully non-linear elliptic equations on compact Hermitian manifolds}, arXiv:1501.02762v3.

\bibitem{GTW} Sz\'ekelyhidi, G., Tosatti, V. and Weinkove, B., {\it Gauduchon metrics with prescribed volume form}, preprint, arXiv:1503.04491.

\bibitem{Wang} Wang, X.J., {\it The k-Hessian equation}, Lect. Not. Math. 1977 (2009).


\bibitem{DZhang} Zhang, D.K., {\it Hessian equations on closed Hermitian manifolds}, arXiv:1501.03553.

\bibitem{Zhang} Zhang, X.W., {\it A priori estimates for complex Monge-Amp\`ere equation on Hermitian manifolds}, Int. Math. Rs. Not. 19 (2010), 3814-3836.

}
\end{thebibliography}
\end{document}